\documentclass[10pt,a4paper]{article}%
\usepackage{amsmath}
\usepackage{amsfonts}
\usepackage[colorlinks=true, pdfstartview=FitV, linkcolor=blue, citecolor=red,
urlcolor=blue]{hyperref}
\usepackage{color}
\usepackage{amssymb}
\usepackage{graphicx}%

\newtheorem{theorem}{Theorem}

\newtheorem{corollary}[theorem]{Corollary}

\newtheorem{lemma}[theorem]{Lemma}

\newtheorem{proposition}[theorem]{Proposition}
\newtheorem{remark}[theorem]{Remark}

\newenvironment{proof}[1][Proof]{\noindent\textbf{#1.} }{\ \rule{0.5em}{0.5em}}

\begin{document}

\title{On the Stieltjes constants with respect to harmonic zeta functions}
\author{Levent Karg\i n, Ayhan Dil, Mehmet Cenkci and M\"{u}m\"{u}n Can\\{\small {Department of Mathematics, Akdeniz University, TR-07058 Antalya,
Turkey}}}
\maketitle

\begin{abstract}
The aim of this paper is to investigate harmonic Stieltjes
constants occurring in the Laurent expansions of the function
\[
\zeta_{H}\left(  s,a\right)  =\sum_{n=0}^{\infty}\frac{1}{\left(  n+a\right)
^{s}}\sum_{k=0}^{n}\frac{1}{k+a},\text{ }\operatorname{Re}\left(  s\right)
>1,
\]
which we call\textit{ harmonic Hurwitz zeta function}. In particular
evaluation formulas for the harmonic Stieltjes constants $\gamma_{H}\left(
m,1/2\right)  $ and $\gamma_{H}\left(  m,1\right)  $ are presented.


\end{abstract}

\section{Introduction}

The Hurwitz zeta function $\zeta\left(  s,a\right)  $ defined by the series%
\[
\zeta\left(  s,a\right)  =\sum_{n=0}^{\infty}\frac{1}{\left(  n+a\right)
^{s}}%
\]
with $\operatorname{Re}\left(  s\right)  >1$ and $a\in\mathbb{C}%
\setminus\left\{  0,-1,-2,\ldots\right\}  $ has an analytic continuation to
the whole complex $s$ plane with the exception at the simple pole $s=1$ with
residue $1$. Near $s=1$, $\zeta\left(  s,a\right)  $ has the Laurent expansion%
\begin{equation}
\zeta\left(  s,a\right)  =\frac{1}{s-1}+\sum_{m=0}^{\infty}\left(  -1\right)
^{m}\frac{\gamma\left(  m,a\right)  }{m!}\left(  s-1\right)  ^{m}, \label{ZL}%
\end{equation}
where the coefficients $\gamma\left(  m,a\right)  $ are called
\textit{generalized} \textit{Stieltjes constants}. In the case $a=1$,
$\zeta\left(  s,1\right)  =\zeta\left(  s\right)  $ is the Riemann zeta
function and the coefficients $\gamma\left(  m,1\right)  =\gamma\left(
m\right)  $ are the \textit{Stieltjes constants} (see for example
\cite{BCh,LiT}).\ The constant $\gamma\left(  0\right)  $ is the famous
\textit{Euler-Mascheroni constant}: $\gamma=\gamma\left(  0\right)
=0.577\,215\,664\,9\ldots$. It is shown that $\gamma\left(  m,a\right)  $ can
be alternatively expressed by the limit%
\begin{equation}
\gamma\left(  m,a\right)  =\lim_{N\rightarrow\infty}\left(  \sum_{n=1}%
^{N}\frac{\log^{m}\left(  n+a\right)  }{n+a}-\frac{\log^{m+1}\left(
N+a\right)  }{m+1}\right)  \label{SC}%
\end{equation}
(see \cite{Be}). There is comprehensive literature on deriving series and
integral representations for the Stieltjes constants and their extensions (see
for example \cite{Bl1,Bl2,CDKCG,Ch,Cof1,Cof2,C,FB,HK}). These representations
usually allow a more accurate estimation of mentioned constants (see
\cite{Ad,Be,Bl2} for more details).

The Dirichlet series associated with the harmonic numbers $H_{n}=\sum
_{k=1}^{n}k^{-1}$, so called the \textit{harmonic zeta function}, is defined
by
\[
\zeta_{H}\left(  s\right)  =\sum_{k=1}^{\infty}\frac{H_{k}}{k^{s}},\text{
}\operatorname{Re}\left(  s\right)  >1,
\]
and subject to many studies. Euler \cite[pp. 217--264]{E} gave a closed form
formula for $\zeta_{H}\left(  s\right)  $ in terms of the Riemann zeta values
for $s\in\mathbb{N}\backslash\left\{  1\right\}  $. Apostol and Vu \cite{AV},
and Matsuoka \cite{Ma} have shown that the function $\zeta_{H}$ can be
continued a meromorphic function with a double pole at $s=1$, and simple poles
at $s=0$ and $s=1-2j,$ $j\in\mathbb{N}$. In \cite{BGP1}, Boyadzhiev et al.
considered the Laurent expansion of the form\textbf{ }%
\begin{equation}
\zeta_{H}\left(  s\right)  =\frac{a_{-1}}{s-b}+a_{0}+O\left(  s-b\right)  ,
\label{zhL}%
\end{equation}
where $b=0$ and $b=1-2j$, $j\in\mathbb{N}$, and gave explicit description for
the coefficient $a_{0}$. Employing the Ramanujan summation method,
Candelpergher and Coppo \cite{CC} have recorded that the \textit{harmonic
Stieltjes constants} $\gamma_{H}\left(  m\right)  $, which occur in the
Laurent expansion
\[
\zeta_{H}\left(  s\right)  =\frac{1}{\left(  s-1\right)  ^{2}}+\frac{\gamma
}{s-1}+\sum_{m=0}^{\infty}\left(  -1\right)  ^{m}\frac{\gamma_{H}\left(
m\right)  }{m!}\left(  s-1\right)  ^{m},\text{ }0<\left\vert s-1\right\vert
<1,
\]
can be expressed as%
\[
\gamma_{H}\left(  m\right)  =\lim_{x\rightarrow\infty}\left(  \sum_{n\leq
x}\frac{H_{n}\ln^{m}n}{n}-\frac{\ln^{m+2}x}{m+2}-\gamma\frac{\ln^{m+1}x}%
{m+1}\right)
\]
(see also \cite[Theorem 2 with $r=1$]{CDK}). Besides, in \cite{CC}, the
constant $\gamma_{H}\left(  0\right)  $ was expressed explicitly and the
coefficient $a_{0}$ in (\ref{zhL}) was rediscovered.  It worths to mention
that, in \cite{Al}, certain real numbers and log-sine integrals have
been shown as a combination of special values of the harmonic zeta function
and the Riemann zeta function, which provided strong approximations for them.

Recently, Alzer and Choi \cite{AC} have introduced four types of parametric
Euler sums, namely%
\begin{align*}
S_{z,s}^{++}\left(  a,b\right)   &  =\sum_{n=1}^{\infty}\frac{\mathcal{H}%
_{n}^{\left(  z\right)  }\left(  a\right)  }{\left(  n+b\right)  ^{s}},\text{
}S_{z,s}^{+-}\left(  a,b\right)  =\sum_{n=1}^{\infty}\left(  -1\right)
^{n+1}\frac{\mathcal{H}_{n}^{\left(  z\right)  }\left(  a\right)  }{\left(
n+b\right)  ^{s}},\\
S_{z,s}^{-+}\left(  a,b\right)   &  =\sum_{n=1}^{\infty}\frac{\mathcal{A}%
_{n}^{\left(  z\right)  }\left(  a\right)  }{\left(  n+b\right)  ^{s}},\text{
}S_{z,s}^{--}\left(  a,b\right)  =\sum_{n=1}^{\infty}\left(  -1\right)
^{n+1}\frac{\mathcal{A}_{n}^{\left(  z\right)  }\left(  a\right)  }{\left(
n+b\right)  ^{s}},
\end{align*}
where%
\[
\mathcal{H}_{n}^{\left(  z\right)  }\left(  a\right)  =\sum_{k=1}^{n}\frac
{1}{\left(  k+a\right)  ^{z}}\text{ and }\mathcal{A}_{n}^{\left(  z\right)
}\left(  a\right)  =\sum_{k=1}^{n}\frac{\left(  -1\right)  ^{k-1}}{\left(
k+a\right)  ^{z}},
\]
$a,b\in\mathbb{C}\backslash\left\{  -1,-2,-3,\ldots\right\}  $ and $s,$
$z\in\mathbb{C}$ are adjusted so that the involved defining series converge.
They have investigated analytic continuations via the Euler-Maclaurin and
Euler-Boole summation formulas, and given shuffle relations. We note that
$S_{1,s}^{++}\left(  0,0\right)  $ is the harmonic zeta function $\zeta
_{H}\left(  s\right)  $, and refer to \cite{BoBB,EL,LiQ,QSL,SoC} for studies
on several types of Euler sums with parameters. In particular, analytic
continuations of $S_{1,s}^{+-}\left(  0,0\right)  $, $S_{1,s}^{--}\left(
0,0\right)  $ and $S_{1,s}^{-+}\left(  0,0\right)  $ have been investigated by
Boyadzhiev et al. in \cite{BGP2}.

Let $H_{n}\left(  a\right)  $ denote the sum
\[
H_{n}\left(  a\right)  =\sum_{k=0}^{n}\frac{1}{k+a}.
\]
It is obvious that $\mathcal{H}_{n}^{\left(  1\right)  }\left(  a-1\right)
=H_{n-1}\left(  a\right)  $ and $H_{n-1}\left(  1\right)  =H_{n}$. Throughout
this study, an empty sum is regarded to be zero.

In this study, we consider the function%
\begin{equation}
\zeta_{H}\left(  s,a\right)  =\sum_{n=0}^{\infty}\frac{H_{n}\left(  a\right)
}{\left(  n+a\right)  ^{s}}, \label{hhz}%
\end{equation}
where $a\in\mathbb{C}\backslash\left\{  0,-1,-2,-3,\ldots\right\}  $ and
$\operatorname{Re}(s)>1.$ It is evident that $\zeta_{H}\left(  s,a\right)  $
is the function $S_{1,s}^{++}\left(  a-1,a-1\right)  $. Since the series
(\ref{hhz}) is a Hurwitz type Dirichlet series and $\zeta_{H}\left(
s,1\right)  =\zeta_{H}\left(  s\right)  $ is the harmonic zeta function, it is
logical to call $\zeta_{H}\left(  s,a\right)  $ as the \textit{harmonic
Hurwitz zeta function}. It is known that the harmonic Hurwitz zeta function
has a second-order pole at $s=1$ and simple poles at $s=0$ and $s=1-2n,$
$n\in\mathbb{N}$ (cf. \cite[Remark 2.1]{AC}). We first prove that%
\begin{align*}
 &\gamma_{H}\left(  m,a\right) \\
 &=\lim_{x\rightarrow\infty}\left(  \sum_{n\leq x}\frac{\log^{m}\left(
n+a\right)  }{n+a}H_{n}\left(  a\right)  -\frac{\log^{m+2}\left(  x+a\right)
}{m+2}-\gamma\left(  0,a\right)  \frac{\log^{m+1}\left(  x+a\right)  }%
{m+1}\right)  ,
\end{align*}
by modifying the method of Briggs and Buschman \cite{BB}. Here the constants
$\gamma_{H}\left(  m,a\right)  $ (\textit{generalized} \textit{harmonic}
\textit{Stieltjes constants}) are defined by
\[
\zeta_{H}\left(  s,a\right)  =\frac{a_{-2}}{\left(  s-1\right)  ^{2}}%
+\frac{a_{-1}}{s-1}+\sum_{m=0}^{\infty}\left(  -1\right)  ^{m}\frac{\gamma
_{H}\left(  m,a\right)  }{m!}\left(  s-1\right)  ^{m},\text{ }0<\left\vert
s-1\right\vert <1.
\]
Secondly, we give a contour integral representation for $\zeta_{H}\left(
s,a\right)  $. Utilizing this representation we find the Laurent expansion of
$\zeta_{H}\left(  s,a\right)  $ at the simple poles $s=k$, where $k=0$ or
$k=1-2n,$ $n\in\mathbb{N}$. We finally give the evaluation formula
(\ref{ev11})\textbf{ }for $\gamma_{H}\left(  m\right)  =\gamma_{H}\left(
m,1\right)  $ in terms of $\gamma\left(  m\right)  $, $\Gamma^{\left(
m\right)  }\left(  1\right)  $, zeta values and certain integrals. As a
demonstration of this formula, we present the first few values of $\gamma
_{H}\left(  m\right)  $ as follows.
\begin{align*}
\gamma_{H}\left(  0\right)   &  =0.98905599532797\ldots & \gamma_{H}\left(
6\right)   &  =359.53098575125632\ldots\\
\gamma_{H}\left(  1\right)   &  =0.40076122995742\ldots & \gamma_{H}\left(
7\right)   &  =2518.3592429515589\ldots\\
\gamma_{H}\left(  2\right)   &  =0.97130466096680\ldots & \gamma_{H}\left(
8\right)   &  =20153.438413223697\ldots\\
\gamma_{H}\left(  3\right)   &  =2.97258466333065\ldots & \gamma_{H}\left(
9\right)   &  =181410.47173964166\ldots\\
\gamma_{H}\left(  4\right)   &  =11.9405013905470\ldots & \gamma_{H}\left(
10\right)   &  =1814252.3513222240\ldots\\
\gamma_{H}\left(  5\right)   &  =59.8445659187868\ldots & \gamma_{H}\left(
11\right)   &  =19957587.910759595\ldots
\end{align*}
We would like to emphesize that an evaluation formula for $\gamma_{H}\left(
m\right)  $ with $m\geq1$ had not been given in \cite{CDK} and \cite{CC} nor,
as far as we know, in any study. Moreover, the formula (\ref{ev11}) provides a
contrubition to the evaluation of the hyperharmonic Stieltjes constants given
in \cite[Theorem 5]{CDK}.

\section{Laurent expansion in a neighborhood of $s=1$}

In this section we prove the following theorem, which is one of the main
results of this paper.

\begin{theorem}
\label{mteo2}Let $0<a\leq1$. The harmonic Hurwitz zeta function has the
following Laurent expansion in the annulus $0<\left\vert s-1\right\vert <1$%
\begin{equation}
\zeta_{H}\left(  s,a\right)  =\frac{1}{\left(  s-1\right)  ^{2}}+\frac
{\gamma\left(  0,a\right)  }{s-1}+\sum_{m=0}^{\infty}\left(  -1\right)
^{m}\frac{\gamma_{H}\left(  m,a\right)  }{m!}\left(  s-1\right)  ^{m},
\label{0}%
\end{equation}
where the coefficients $\gamma_{H}\left(  m,a\right)  $ are determined by
\begin{align*}
 &\gamma_{H}\left(  m,a\right)  \\
& =\lim_{x\rightarrow\infty}\left(  \sum_{n\leq x}\frac{H_{n}\left(
a\right)  \log^{m}\left(  n+a\right)  }{n+a}-\frac{\log^{m+2}\left(
x+a\right)  }{m+2}-\gamma\left(  0,a\right)  \frac{\log^{m+1}\left(
x+a\right)  }{m+1}\right)  .
\end{align*}

\end{theorem}

\begin{remark}
The definition (\ref{hhz}) of $\zeta_{H}\left(  s,a\right)  $ makes it obvious
that $\zeta_{H}\left(  s,a+1\right)  =\zeta_{H}\left(  s,a\right)  -\frac
{1}{a}\zeta\left(  s,a\right)  $, from which it follows that constants
$\gamma_{H}\left(  m,a\right)  $ satisfy the recurrence relation
\[
\gamma_{H}\left(  m,a+1\right)  =\gamma_{H}\left(  m,a\right)  -\frac{1}%
{a}\gamma\left(  m,a\right)  ,\text{ }a\not =0,-1,-2,\ldots,\text{
}m=0,1,2,\ldots
\]
Moreover, for a positive integer $n$, we can write%
\[
\gamma_{H}\left(  m,n+1\right)  =\gamma_{H}\left(  m\right)  -\sum_{k=1}%
^{n}\frac{\gamma\left(  m,k\right)  }{k}.
\]

\end{remark}

The basis of our approach depends on the equation
\begin{align}
\zeta_{H}\left(  s,a\right)   &  =\frac{\zeta\left(  s,a\right)  }{s-1}%
+\frac{1}{2}\zeta\left(  s+1,a\right)  +\sum_{j=1}^{q}\frac{B_{2j}\zeta\left(
s+2j,a\right)  }{\left(  2j\right)  !}s\left(  s+1\right)  \cdots\left(
s+2j-2\right) \nonumber\\
&  -\frac{s\left(  s+1\right)  \cdots\left(  s+2q\right)  }{\left(
2q+1\right)  !}\sum_{k=0}^{\infty}\frac{1}{k+a}%
%TCIMACRO{\dint _{k}^{\infty}}%
%BeginExpansion
{\displaystyle\int_{k}^{\infty}}
%EndExpansion
\frac{\overline{B}_{2q+1}\left(  t\right)  }{\left(  t+a\right)  ^{s+2q+1}}dt,
\label{hzeta}%
\end{align}
which is equation (2.16) of \cite{AC} with $p=1$. Here $\overline{B}%
_{n}\left(  x\right)  =B_{n}\left(  x-\left\lfloor x\right\rfloor \right)  $,
$\left\lfloor x\right\rfloor $ is the greatest integer not exceeding $x$ and
$B_{n}\left(  x\right)  $ are the Bernoulli polynomials defined by the
generating function%
\begin{equation}
\frac{ze^{xz}}{e^{z}-1}=\sum_{n=0}^{\infty}\frac{B_{n}\left(  x\right)  }%
{n!}z^{n}. \label{Bn}%
\end{equation}
The case $x=0$ are called Bernoulli numbers and denoted by $B_{n}=B_{n}\left(
0\right)  $ with $B_{0}=1$, $B_{1}=-1/2$, $B_{2j+1}=0$ for $j\geq1$. Now,
(\ref{hzeta}) makes it obvious that the Laurent series of the harmonic Hurwitz
zeta function can be written as
\[
\zeta_{H}\left(  s,a\right)  =\frac{1}{\left(  s-1\right)  ^{2}}+\frac
{\gamma\left(  0,a\right)  }{s-1}+\sum_{m=0}^{\infty}a_{m}\left(  s-1\right)
^{m},\text{ }0<|s-1|<1.
\]
Therefore the proof of Theorem \ref{mteo2} will follow once we show that
\[
a_{m}   =\frac{\left(  -1\right)  ^{m}}{m!}\lim_{x\rightarrow\infty}\left(  \sum_{n\leq x}\frac{H_{n}\left(
a\right)  \log^{m}\left(  n+a\right)  }{n+a}-\frac{\log^{m+2}\left(
x+a\right)  }{m+2}-\gamma\left(  0,a\right)  \frac{\log^{m+1}\left(
x+a\right)  }{m+1}\right)  .
\]
We reach at this result by utilizing some auxiliary statements.

\begin{lemma}
[Abel summation formula](see \cite[Theorem 4.2]{Ap}) \label{lem1}If $\left\{
b_{1},b_{2},b_{3},\ldots\right\}  $ is a sequence of complex numbers and
$v\left(  x\right)  $ is a function with a continuous derivative for
$x>\alpha$, then%
\[
\sum_{n\leq x}b_{n}v\left(  n\right)  =\left(  \sum_{n\leq x}b_{n}\right)
v\left(  x\right)  -\int\limits_{\alpha}^{x}\left(  \sum_{n\leq t}%
b_{n}\right)  v^{\prime}\left(  t\right)  dt.
\]

\end{lemma}

Setting $b_{n}=H_{n}\left(  a\right)  $ and $v\left(  x\right)  =\left(
x+a\right)  ^{-s}$ in Lemma \ref{lem1} gives the following result.

\begin{lemma}
\label{lem2}For $\operatorname{Re}\left(  s\right)  >1$ we have%
\[
\zeta_{H}\left(  s,a\right)  =s\int\limits_{1-a}^{\infty}\left(  t+a\right)
^{-s-1}\left(  \sum_{n\leq t}H_{n}\left(  a\right)  \right)  dt.
\]

\end{lemma}

Now as $k\rightarrow\infty$, we have the well-known identity%
\[
H_{k}\left(  a\right)  =\log\left(  k+a\right)  +\gamma\left(  0,a\right)
+O\left(  \frac{1}{k+a}\right)  ,
\]
so the expression%
\[
\sum_{n=0}^{k}H_{n}\left(  a\right)  =\sum_{j=0}^{k}\frac{1}{j+a}\sum
_{n=j}^{k}1=\left(  k+1+a\right)  H_{k}\left(  a\right)  -\left(  k+1\right)
\]
yields to%
\[
\sum_{n=0}^{k}H_{n}\left(  a\right)  =\left(  k+a\right)  \log\left(
k+a\right)  +\left(  k+a\right)  \left(  \gamma\left(  0,a\right)  -1\right)
+O\left(  \log\left(  k+a\right)  \right)  .
\]
Let%
\begin{align}
E\left(  x,a\right)   &  =\sum_{n\leq x}H_{n}\left(  a\right)  -\left(
x+a\right)  \log\left(  x+a\right)  -\left(  x+a\right)  \left(  \gamma\left(
0,a\right)  -1\right) \nonumber\\
&  =O\left(  \log\left(  x+a\right)  \right)  . \label{5}%
\end{align}

\begin{lemma}
\label{lem3}For $\operatorname{Re}\left(  s\right)  >0$ we have%
\[
f\left(  s,a\right)  :=s\int\limits_{1-a}^{\infty}\left(  t+a\right)
^{-s-1}E\left(  t,a\right)  dt=1-\gamma\left(  0,a\right)  +\sum_{n=0}%
^{\infty}a_{n}\left(  s-1\right)  ^{n}.
\]

\end{lemma}

\begin{proof}
For $\operatorname{Re}\left(  s\right)  >0$ the integral is an analytic
function. Moreover, for $\operatorname{Re}\left(  s\right)  >1$, it follows
from (\ref{5}) and Lemma \ref{lem2} that%
\begin{align*}
  &s\int\limits_{1-a}^{\infty}\left(  x+a\right)  ^{-s-1}E\left(  x,a\right)
dx\\
&  =\zeta_{H}\left(  s,a\right)  -s\int\limits_{1-a}^{\infty}\left(
x+a\right)  ^{-s}\log\left(  x+a\right)  dx-s\left(  \gamma\left(  0,a\right)
-1\right)  \int\limits_{1-a}^{\infty}\left(  x+a\right)  ^{-s}dx\\
&  =1-\gamma\left(  0,a\right)  +\sum_{n=0}^{\infty}a_{n}\left(  s-1\right)
^{n},
\end{align*}
which is the desired conclusion.
\end{proof}

Next result plays a crucial role in proving Theorem \ref{mteo2}.

\begin{theorem}
\label{teo1}Let $m$ be a non-negative integer and let $u<0$. Then,%
\begin{align*}
\sum_{n\leq x}\frac{H_{n}\left(  a\right)  \log^{m}\left(  n+a\right)
}{\left(  n+a\right)  ^{-u}}  &  =\int\limits_{1-a}^{x}\frac{\log^{m+1}\left(
t+a\right)  +\gamma\left(  0,a\right)  \log^{m}\left(  t+a\right)  }{\left(
t+a\right)  ^{-u}}dt\\
&  +\left(  -1\right)  ^{m}f^{\left(  m\right)  }\left(  -u,a\right)
+o\left(  1\right)  +\left\{
\begin{array}
[c]{lc}%
\gamma\left(  0,a\right)  -1, & m=0,\\
0, & m=1,
\end{array}
\right.
\end{align*}
where $f\left(  s,a\right)  $ is the function introduced in Lemma \ref{lem3}.
\end{theorem}

\begin{proof}
Set $v\left(  x\right)  =\left(  x+a\right)  ^{u}\log^{m}\left(  x+a\right)  $
and $b_{n}=H_{n}\left(  a\right)  $ in Lemma \ref{lem1}. Then%
\begin{align}
\sum_{n\leq x}\frac{H_{n}\left(  a\right)  \log^{m}\left(  n+a\right)
}{\left(  n+a\right)  ^{-u}}  &  =\left(  x+a\right)  ^{u}\log^{m}\left(
x+a\right)  \sum_{n\leq x}H_{n}\left(  a\right) \nonumber\\
&  -\int\limits_{1-a}^{x}\sum_{n\leq t}H_{n}\left(  a\right)  \frac{d}%
{dt}\left(  \left(  t+a\right)  ^{u}\log^{m}\left(  t+a\right)  \right)  dt.
\label{9}%
\end{align}
Making use of the equality
\[
\frac{d}{dt}\left(  \left(  t+a\right)  ^{u}\log^{m}\left(  t+a\right)
\right)  =\frac{d^{m}}{du^{m}}\left(  u\left(  t+a\right)  ^{u-1}\right)
\]
enables us to write the RHS of (\ref{9}) as%
\[
L_{1}+\left(  \gamma\left(  0,a\right)  -1\right)  L_{2}-\frac{d^{m}}{du^{m}%
}u\int\limits_{1-a}^{x}E\left(  t,a\right)  \left(  t+a\right)  ^{u-1}%
dt+O\left(  x^{u}\log^{m+1}x\right)  ,
\]
where%
\begin{align*}
L_{1}  &  =\left(  x+a\right)  ^{u+1}\log^{m+1}\left(  x+a\right)
-\frac{d^{m}}{du^{m}}u\int\limits_{1-a}^{x}\left(  t+a\right)  ^{u}\log\left(
t+a\right)  dt,\\
L_{2}  &  =\left(  x+a\right)  ^{u+1}\log^{m}\left(  x+a\right)  -\frac{d^{m}%
}{du^{m}}u\int\limits_{1-a}^{x}\left(  t+a\right)  ^{u}dt.
\end{align*}
It is a simple matter to see that\textbf{ }%
\[
L_{2}=\frac{d^{m}}{du^{m}}\int\limits_{1-a}^{x}\left(  t+a\right)
^{u}dt+\frac{d^{m}}{du^{m}}1=\int\limits_{1-a}^{x}\left(  t+a\right)  ^{u}%
\log^{m}\left(  t+a\right)  dt+\left\{
\begin{array}
[c]{cc}%
1, & m=0\\
0, & m=1
\end{array}
\right.
\]
and
\begin{align*}
L_{1}  &  =\frac{d^{m}}{du^{m}}\left(  \int\limits_{1-a}^{x}\left(
t+a\right)  ^{u}\log\left(  t+a\right)  dt+\int\limits_{1-a}^{x}\left(
t+a\right)  ^{u}dt\right) \\
&  =\int\limits_{1-a}^{x}\left(  t+a\right)  ^{u}\log^{m+1}\left(  t+a\right)
dt+\int\limits_{1-a}^{x}\left(  t+a\right)  ^{u}\log^{m}\left(  t+a\right)
dt.
\end{align*}
If we put these together, we get%
\begin{align*}
&\sum_{n\leq x}\frac{H_{n}\left(  a\right)  \log^{m}\left(  n+a\right)
}{\left(  n+a\right)  ^{-u}}\\
&  =\int\limits_{1-a}^{x}\frac{\log^{m+1}\left(
t+a\right)  +\gamma\left(  0,a\right)  \log^{m}\left(  t+a\right)  }{\left(
t+a\right)  ^{-u}}dt\\
&  +\left(  -1\right)  ^{m}f^{\left(  m\right)  }\left(  -u,a\right)
+O\left(  x^{u}\log^{m+1}x\right)  +\left\{
\begin{array}
[c]{cc}%
\gamma\left(  0,a\right)  -1, & m=0,\\
0, & m=1,
\end{array}
\right.
\end{align*}
and the proof is complete.
\end{proof}

We are now ready to prove the main result of this section.

\begin{proof}
[Proof of Theorem \ref{mteo2}]From Lemma \ref{lem3} we have
\begin{equation}
f\left(  1,a\right)  =1-\gamma\left(  0,a\right)  +a_{0}\text{ and }%
a_{m}=\frac{f^{\left(  m\right)  }\left(  1,a\right)  }{m!}\text{ for }m>0.
\label{14}%
\end{equation}
Setting $u=-1$ in Theorem \ref{teo1} yields to
\begin{align}
\sum_{n\leq x}\frac{H_{n}\left(  a\right)  }{n+a}\log^{m}\left(  n+a\right)
&  =\int\limits_{1-a}^{x}\frac{\log^{m+1}\left(  t+a\right)  +\gamma\left(
0,a\right)  \log^{m}\left(  t+a\right)  }{t+a}dt\nonumber\\
&  +\left(  -1\right)  ^{m}f^{\left(  m\right)  }\left(  1,a\right)  +o\left(
1\right)  +\left\{
\begin{array}
[c]{cc}%
\hspace{-0.06in}\gamma\left(  0,a\right)  -1, & \hspace{-0.1in}m=0,\\
\hspace{-0.06in}0, & \hspace{-0.1in}m=1.
\end{array}
\right.  \label{15}%
\end{align}
Hence (\ref{14}) and (\ref{15}) lead to%
\begin{align*}
&a_{m}    =\frac{\left(  -1\right)  ^{m}}{m!}\\
&\times\lim_{x\rightarrow\infty}\left(  \sum_{n\leq x}\frac{H_{n}\left(
a\right)  \log^{m}\left(  n+a\right)  }{n+a}-\frac{\log^{m+2}\left(
x+a\right)  }{m+2}-\gamma\left(  0,a\right)  \frac{\log^{m+1}\left(
x+a\right)  }{m+1}\right)  ,
\end{align*}
which is our claim.
\end{proof}

\section{Laurent expansion via contur integral}

The main result of Boyadzhiev et al.'s paper \cite{BGP1} is the formula%
\begin{equation}
\zeta_{H}\left(  s\right)  =\Gamma\left(  1-s\right)  I\left(  s\right)
+\zeta\left(  s+1\right)  -\psi\left(  s\right)  \zeta\left(  s\right)
-\zeta^{\prime}\left(  s\right)  , \label{B1}%
\end{equation}
where $\Gamma\left(  s\right)  $ is the Euler gamma function and
\[
I\left(  s\right)  =\frac{1}{2\pi i}\int_{C}\frac{z^{s-1}e^{z}}{e^{z}-1}%
\log\left(  \frac{e^{z}-1}{z}\right)  dz.
\]
Here $C$ denotes the Hankel contour which starts from $-\infty$ along the
lower side of the negative real axis, encircles the origin once in the
positive (counter-clockwise) direction and then returns to $-\infty$ along the
upper side of the negative real axis. The loop $C$ consists of the parts
$C=C_{-}\cup C_{+}\cup C_{\varepsilon}$, where $C_{\varepsilon}$ is a
positively-oriented circle of radius $\varepsilon$ about the origin, and
$C_{-}$ and $C_{+}$ are the lower and upper edges of a cut in the complex
$z$-plane along the negative real axis.

In this section we first give a contour integral representation for $\zeta
_{H}\left(  s,a\right)  $. This is motivated by \cite{BGP1}. We next study the
special values of $\zeta_{H}\left(  -2k,a\right)  $ for positive integers $k$
and Laurent series expansions at the simple poles with the help of this
integral representation.

In what follows, we denote the Hankel contour by $C=C_{-}\cup C_{+}\cup
C_{\varepsilon}$ and use the parameterizations $z=xe^{-\pi i}$ on $C_{-}$,
$z=xe^{\pi i}$ on $C_{+}$, and $z=\varepsilon e^{i\theta}$ on $C_{\varepsilon
}$.

Before giving the integral representation we would like to point out that
$H_{n}\left(  a\right)  $ has the the following generating function%
\[
\sum_{n=0}^{\infty}H_{n}\left(  a\right)  x^{n}=\frac{\Phi\left(
x,1;a\right)  }{1-x},\text{ }\left\vert x\right\vert <1,
\]
where%
\[
\Phi\left(  x,s;a\right)  =\sum_{k=0}^{\infty}\frac{x^{k}}{\left(  k+a\right)
^{s}}%
\]
is the Lerch transcendent (see for example \cite[Sec. 2.5]{SC}).

\begin{theorem}
\label{mteo1}Let $a\in\mathbb{C}\setminus\left\{  0,-1,-2,\ldots\right\}  $.
For $\operatorname{Re}\left(  s\right)  >1$ we have%
\begin{equation}
\zeta_{H}\left(  s,a\right)  =\Gamma\left(  1-s\right)  I\left(  s,a\right)
+\zeta_{H}\left(  s\right)  , \label{2}%
\end{equation}
where
\[
I\left(  s,a\right)  =\frac{1}{2\pi i}\int_{C}\frac{z^{s-1}}{1-e^{z}}\left\{
e^{za}\Phi\left(  e^{z},1;a\right)  +\log\left(  1-e^{z}\right)  \right\}  dz
\]
and $C$ is the aforementioned Hankel contour.\textbf{\ }
\end{theorem}

\begin{proof}
Considering%
\[
\Gamma\left(  s\right)  =\int_{0}^{\infty}t^{s-1}e^{-t}dt
\]
we can write the harmonic Hurwitz zeta function as%
\begin{align}
\zeta_{H}\left(  s,a\right)   &  =\sum_{n=0}^{\infty}\frac{H_{n}\left(
a\right)  }{\left(  n+a\right)  ^{s}}=\frac{1}{\Gamma\left(  s\right)  }%
\int_{0}^{\infty}x^{s-1}e^{-ax}\sum_{n=0}^{\infty}H_{n}\left(  a\right)
e^{-nx}dx\nonumber\\
&  =\frac{1}{\Gamma\left(  s\right)  }\int_{0}^{\infty}\frac{x^{s-1}e^{-xa}%
}{1-e^{-x}}\Phi\left(  e^{-x},1;a\right)  dx \label{10}%
\end{align}
for $\operatorname{Re}\left(  s\right)  >1$.

We now focus on the integral $I\left(  s,a\right)  .$ Since $z=xe^{-\pi i}$ on
$C_{-}$, $z=xe^{\pi i}$ on $C_{+}$ and $z=\varepsilon e^{i\theta}$ on
$C_{\varepsilon}$ we find that%
\begin{align*}
I\left(  s,a\right)   &  =\frac{1}{2\pi i}\left(  \int_{C_{-}}+\int_{C_{+}%
}+\int_{C_{\varepsilon}}\right)  \frac{z^{s-1}}{1-e^{z}}\left\{  e^{za}%
\Phi\left(  e^{z},1;a\right)  +\log\left(  1-e^{z}\right)  \right\}  dz\\
&  =\frac{\sin\left(  s\pi\right)  }{\pi}\int_{\varepsilon}^{\infty}%
\frac{x^{s-1}}{1-e^{-x}}\left(  e^{-xa}\Phi\left(  e^{-x},1;a\right)
+\log\left(  1-e^{-x}\right)  \right)  dx\\
&  +\frac{1}{2\pi i}\int_{C_{\varepsilon}}\frac{z^{s-1}}{1-e^{z}}\left\{
e^{za}\Phi\left(  e^{z},1;a\right)  +\log\left(  1-e^{z}\right)  \right\}  dz.
\end{align*}
To show that%
\[
\left\vert \frac{1}{2\pi i}\int_{C_{\varepsilon}}\frac{z^{s-1}}{1-e^{z}%
}\left\{  e^{za}\Phi\left(  e^{z},1;a\right)  +\log\left(  1-e^{z}\right)
\right\}  dz\right\vert \longrightarrow0\text{ as }\varepsilon\rightarrow0
\]
we utilize the expression (see for example \cite[Sec. 2.5]{SC})%
\begin{equation}
\Phi\left(  e^{z},1;a\right)  =e^{-za}\left\{  \sum_{n=1}^{\infty}\frac
{\zeta\left(  1-n,a\right)  }{n!}z^{n}+\psi\left(  1\right)  -\psi\left(
a\right)  -\log\left(  -z\right)  \right\}  , \label{8}%
\end{equation}
where $\psi\left(  z\right)  =\Gamma^{\prime}\left(  z\right)  /\Gamma\left(
z\right)  $ is the digamma function. We then have%
\begin{align*}
&  \frac{1}{2\pi i}\int_{C_{\varepsilon}}\frac{z^{s-1}}{1-e^{z}}\left\{
e^{za}\Phi\left(  e^{z},1;a\right)  +\log\left(  1-e^{z}\right)  \right\}
dz\\
&  =\sum_{n=1}^{\infty}\frac{\zeta\left(  1-n,a\right)  }{n!}\frac{1}{2\pi
i}\int_{C_{\varepsilon}}\frac{z^{s-1}}{1-e^{z}}z^{n}dz-\frac{1}{2\pi i}%
\int_{C_{\varepsilon}}\frac{z^{s-1}}{1-e^{z}}\log\left(  \frac{-z}{1-e^{z}%
}\right)  dz\\
&  +\left(  \psi\left(  1\right)  -\psi\left(  a\right)  \right)  \frac
{1}{2\pi i}\int_{C_{\varepsilon}}\frac{z^{s-1}}{1-e^{z}}dz.
\end{align*}
For $\operatorname{Re}\left(  s\right)  >1$ the functions $\frac{z^{s-1}%
}{1-e^{z}}z^{n}$ and $\frac{z^{s-1}}{1-e^{z}}\log\left(  \frac{-z}{1-e^{z}%
}\right)  $ are holomorphic in a neighborhood of zero. So the first two
integrals vanish. Besides it is known that the third integral tends to zero as
$\varepsilon\rightarrow0$ (see \cite[Chapter 2]{SC}). Accordingly, we reach at%
\begin{equation}
I\left(  s,a\right)  =\frac{\sin\left(  s\pi\right)  }{\pi}\int_{0}^{\infty
}\frac{x^{s-1}}{1-e^{-x}}\left(  e^{-xa}\Phi\left(  e^{-x},1;a\right)
+\log\left(  1-e^{-x}\right)  \right)  dx. \label{Isa}%
\end{equation}
The proof now follows from (\ref{10}) with the use of $\Gamma(s)\Gamma
(1-s)=\frac{\pi}{\sin(s\pi)},$ $e^{-x}\Phi(e^{-x},\allowbreak1;\allowbreak
1)=-\log\left(  1-e^{-x}\right)  $ and $\zeta_{H}\left(  s,1\right)
=\zeta_{H}\left(  s\right)  .$
\end{proof}

When we consider Theorem \ref{mteo1} for non-positive integer values of $s$,
we need the coefficients $A_{n}\left(  x\right)  $ defined by the generating
function%
\begin{equation}
\frac{ze^{xz}}{e^{z}-1}\log\left(  \frac{e^{z}-1}{z}\right)  =\sum
_{n=1}^{\infty}A_{n}\left(  x\right)  z^{n}. \label{An}%
\end{equation}
The case $x=1$ is nothing but Eq. (6) of \cite{BGP1}. Note that since (cf.
\cite{BGP1})
\[
\log\left(  \frac{e^{z}-1}{z}\right)  =\sum_{n=1}^{\infty}\frac{B_{n}\left(
1\right)  }{n!n}z^{n},
\]
$A_{n}\left(  x\right)  $ can be written as%
\[
A_{n}\left(  x\right)  =\frac{1}{n!}\sum_{k=1}^{n}\binom{n}{k}\frac{1}%
{k}B_{n-k}\left(  x\right)  B_{k}\left(  1\right)  .
\]

\begin{corollary}
For $k\in\mathbb{N},$ we have%
\begin{align*}
\zeta_{H}\left(  -2k,a\right)   & =-\frac{1}{2k+1}\sum_{n=0}^{k}\binom
{2k+1}{2n+1}\zeta\left(  -2n,a\right)  B_{2k-2n}+\frac{1}{2}\zeta\left(  1-2k,a\right) \\
& +\frac{B_{2k}}{2}-\frac
{B_{2k}}{4k}-\left(  2k\right)  !A_{2k+1}\left(  0\right)  .
\end{align*}

\end{corollary}

\begin{proof}
Theorem \ref{mteo1} yields to%
\begin{equation}
\zeta_{H}\left(  -2k,a\right)  =\Gamma\left(  2k+1\right)  I\left(
-2k,a\right)  +\zeta_{H}\left(  -2k\right)  \label{zeta2k}%
\end{equation}
for $k\in\mathbb{N}.$ It is known that $\zeta_{H}\left(  -2k\right)
=\frac{B_{2k}}{2}-\frac{B_{2k}}{4k}$ (see \cite{Ma} or \cite{BGP1}). Using the
residue theorem with (\ref{Bn}), (\ref{8}) and (\ref{An}) the values $I\left(
-2k,a\right)  $ can be evaluated as%
\[
I\left(  -2k,a\right)  =-\sum_{\substack{m+n=2k+1\\n\geq1,m\geq0}}\frac
{\zeta\left(  1-n,a\right)  }{n!}\frac{B_{m}}{m!}-\left(  \psi\left(
1\right)  -\psi\left(  a\right)  \right)  \frac{B_{2k+1}}{\left(  2k+1\right)
!}-A_{2k+1}\left(  0\right)  .
\]
Since $B_{2j+1}=0$ for integers $j\geq1$, we have%
\[
I\left(  -2k,a\right)  =-\sum_{n=0}^{k}\frac{\zeta\left(  -2n,a\right)
}{\left(  2n+1\right)  !}\frac{B_{2k-2n}}{\left(  2k-2n\right)  !}-\frac
{\zeta\left(  1-2k,a\right)  }{\left(  2k\right)  !}B_{1}-A_{2k+1}\left(
0\right)  ,
\]
which gives the assertion of the corollary when substituted in (\ref{zeta2k}).
\end{proof}

We now consider the cases $s=0$ and $s=k=1-2j,$\ $j\in\mathbb{N},$\ in Theorem
\ref{mteo1} to obtain Laurent series expansions of the harmonic Hurwitz zeta
function in the form\
\[
\zeta_{H}\left(  s,a\right)  =\frac{a_{-1}}{s-k}+a_{0}+O\left(  s-k\right)  .
\]

\begin{corollary}
In a neighborhood of zero we have%
\[
\zeta_{H}\left(  s,a\right)  =\frac{1}{2s}+\gamma_{H,0}\left(  0,a\right)
+O\left(  s\right)  ,
\]
where
\[
\gamma_{H,0}\left(  0,a\right)  =-\zeta\left(  0,a\right)  -\frac{1}{2}%
\psi\left(  a\right)
\]

\end{corollary}

\begin{proof}
From Theorem \ref{mteo1} and
\[
\zeta_{H}\left(  s\right)  =\frac{1}{2s}+\frac{1+\gamma\left(  0\right)  }%
{2}+O\left(  s\right)
\]
(see \cite[Corollary 2]{BGP1}) we have%
\[
\zeta_{H}\left(  s,a\right)  =\frac{1}{2s}+\frac{1+\gamma\left(  0\right)
}{2}+I\left(  0,a\right)  +O\left(  s\right)  .
\]
Using the residue theorem with (\ref{Bn}), (\ref{8}) and (\ref{An}) the value
$I\left(  0,a\right)  $ is computed to be
\[
I\left(  0,a\right)  =-\zeta\left(  0,a\right)  -\left(  \psi\left(  1\right)
-\psi\left(  a\right)  \right)  B_{1}-A_{1}\left(  0\right)  .
\]
The proof then is completed.
\end{proof}

\begin{corollary}
In a neighborhood of $s=-1$ we have%
\[
\zeta_{H}\left(  s,a\right)  =-\frac{1}{12\left(  s+1\right)  }+\gamma
_{H,-1}\left(  0,a\right)  +O\left(  s+1\right)  ,
\]
where%
\[
\gamma_{H,-1}\left(  0,a\right)  =-\frac{1}{8}+\frac{\gamma\left(  0\right)
}{6}+\frac{\psi\left(  a\right)  }{12}+\frac{1}{2}\left(  \zeta\left(
0,a\right)  -\zeta\left(  -1,a\right)  \right)  -A_{2}\left(  0\right)  .
\]

\end{corollary}

\begin{proof}
From Theorem \ref{mteo1} and
\[
\zeta_{H}\left(  s\right)  =-\frac{1}{12\left(  s+1\right)  }-\frac{1}%
{8}+\frac{\gamma\left(  0\right)  }{12}+O\left(  s+1\right)
\]
(\cite[Corollary 3]{BGP1}) it is seen that%
\[
\zeta_{H}\left(  s,a\right)  =-\frac{1}{12\left(  s+1\right)  }-\frac{1}%
{8}+\frac{\gamma\left(  0\right)  }{12}+I\left(  -1,a\right)  +O\left(
s+1\right)  .
\]
Using (\ref{Bn}), (\ref{8}) and (\ref{An}), we find that%
\[
I\left(  -1,a\right)  =-\sum_{\substack{m+n=2\\n\geq1,m\geq0}}\frac
{\zeta\left(  1-n,a\right)  }{n!}\frac{B_{m}}{m!}-\left(  \psi\left(
1\right)  -\psi\left(  a\right)  \right)  \frac{B_{2}}{2!}-A_{2}\left(
0\right)  ,
\]
which yields to the desired result.
\end{proof}

Finally we state an expansion in the neighborhood of $s=1-2k$ for integers
$k\geq2$.

\begin{corollary}
Let $k\geq2$ be an integer. In a neighborhood of $s=1-2k$, we have%
\[
\zeta_{H}\left(  s,a\right)  =\frac{\zeta\left(  1-2k\right)  }{s+2k-1}%
+\gamma_{H,1-2k}\left(  0,a\right)  +O\left(  s+2k-1\right)  ,
\]
where
\begin{align*}
\gamma_{H,1-2k}\left(  0,a\right)    & =\frac{1}{2}\zeta\left(  2-2k,a\right)
-\psi\left(  2k\right)  \zeta\left(  1-2k\right)  \\
&-\left(  \psi\left(  1\right)  -\psi\left(  a\right)  \right)  \frac
{B_{2k}}{2k}-\frac{1}{2k}\sum_{n=1}^{k}\dbinom{2k}{2n}\zeta\left(
1-2n,a\right)  B_{2k-2n}.
\end{align*}

\end{corollary}

\begin{proof}
From Theorem \ref{mteo1} and%
\[
\zeta_{H}\left(  s\right)  =\frac{\zeta\left(  1-2k\right)  }{s+2k-1}%
+\Gamma\left(  2k\right)  A_{2k}\left(  1\right)  -\psi\left(  2k\right)
\zeta\left(  1-2k\right)  +O\left(  s+2k-1\right)  ,\text{ }2\leq
k\in\mathbb{N}%
\]
(\cite[Corollary 3]{BGP1}) we have%
\begin{align*}
\zeta_{H}\left(  s,a\right)    &=\frac{\zeta\left(  1-2k\right)  }%
{s+2k-1}+\Gamma\left(  2k\right)  A_{2k}\left(  1\right) \\
& -\psi\left(
2k\right)  \zeta\left(  1-2k\right)  +\Gamma\left(  2k\right)  I\left(  1-2k,a\right)  +O\left(  s+2k-1\right)
.
\end{align*}
Again from residue theorem, (\ref{Bn}), (\ref{8}) and (\ref{An}) we deduce
that%
\begin{align*}
I\left(  1-2k,a\right)     &=-\sum_{n=1}^{k}\frac{\zeta\left(  1-2n,a\right)
}{2n!}\frac{B_{2k-2n}}{\left(  2k-2n\right)  !} \\
&-\frac{\zeta\left(  2-2k,a\right)  }{\left(  2k-1\right)  !}B_{1}-\left(
\psi\left(  1\right)  -\psi\left(  a\right)  \right)  \frac{B_{2k}}{\left(
2k\right)  !}-A_{2k}\left(  0\right)  ,
\end{align*}
from which the proof follows.
\end{proof}

\section{Evaluation formulas for $\gamma_{H}\left(  m,1/2\right)  $ and
$\gamma_{H}\left(  m\right)  $}

According to \cite[Theorems 1 and 3]{BGP2} we have the following
Taylor/Laurent expansions:
\begin{equation}%
\begin{array}
[c]{c}%
\zeta_{\widetilde{H}}\left(  s\right)  =S_{1,s}^{+-}\left(  0,0\right)  =%
%TCIMACRO{\dsum \limits_{n=1}^{\infty}}%
%BeginExpansion
{\displaystyle\sum\limits_{n=1}^{\infty}}
%EndExpansion
\left(  -1\right)  ^{n-1}\dfrac{H_{n}}{n^{s}}=%
%TCIMACRO{\dsum \limits_{m=0}^{\infty}}%
%BeginExpansion
{\displaystyle\sum\limits_{m=0}^{\infty}}
%EndExpansion
\left(  -1\right)  ^{m}\dfrac{\gamma_{\widetilde{H}}\left(  m\right)  }%
{m!}\left(  s-1\right)  ^{m},\medskip\\
\zeta_{H^{^{-}}}\left(  s\right)  =S_{1,s}^{-+}\left(  0,0\right)  =%
%TCIMACRO{\dsum \limits_{n=1}^{\infty}}%
%BeginExpansion
{\displaystyle\sum\limits_{n=1}^{\infty}}
%EndExpansion
\dfrac{H_{n}^{^{-}}}{n^{s}}=\dfrac{\log2}{s-1}+%
%TCIMACRO{\dsum \limits_{m=0}^{\infty}}%
%BeginExpansion
{\displaystyle\sum\limits_{m=0}^{\infty}}
%EndExpansion
\left(  -1\right)  ^{m}\dfrac{\gamma_{H^{^{-}}}\left(  m\right)  }{m!}\left(
s-1\right)  ^{m},
\end{array}
\label{12a}%
\end{equation}
where $H_{n}^{^{-}}=\mathcal{A}_{n}^{\left(  1\right)  }\left(  0\right)  $,
the $n$th skew harmonic number.

In this section, we first evaluate the constants $\gamma_{H}\left(
m,1/2\right)  $ in terms of $\gamma_{H}\left(  m\right)  ,$ $\gamma
_{\widetilde{H}}\left(  k\right)  $ and $\gamma_{H^{^{-}}}\left(  k\right)  $.
We then give an evaluation formula for $\gamma_{H}\left(  m\right)  $ in terms
of $\gamma\left(  m\right)  $, $\Gamma^{\left(  m\right)  }\left(  1\right)
$, zeta values and certain integrals.

\begin{proposition}
For a non-negative integer $m$ we have%
\[
\gamma_{H}\left(  m,\frac{1}{2}\right)  =\gamma_{H}\left(  m\right)
+2\left(  -1\right)  ^{m}\frac{\log^{m+1}2}{m+1}  +2\sum_{k=0}^{m}\binom{m}{k}\left(  -\log2\right)  ^{m-k}\left(
\gamma_{\widetilde{H}}\left(  k\right)  +\gamma_{H^{^{-}}}\left(  k\right)
\right)  .
\]

\end{proposition}

\begin{proof}
From Theorem \ref{mteo1} with (\ref{Isa}) we have%
\begin{align*}
\zeta_{H}\left(  s,1\right)  +\zeta_{H}\left(  s,\frac{1}{2}\right)   &
=2\zeta_{H}\left(  s\right)  +\frac{2}{\Gamma\left(  s\right)  }\int%
_{0}^{\infty}\frac{x^{s-1}}{1-e^{-x}}\log\left(  \frac{1-e^{-x}}{1-e^{-x/2}%
}\right)  dx\\
&  =2\zeta_{H}\left(  s\right)  +\frac{2}{\Gamma\left(  s\right)  }\int%
_{0}^{\infty}\frac{x^{s-1}}{1-e^{-x}}\log\left(  1+e^{-x/2}\right)  dx.
\end{align*}
Employing the partial fractions
\[
\frac{1}{1-t^{2}}=\frac{1}{2}\left(  \frac{1}{1-t}+\frac{1}{1+t}\right)
\]
and considering the generating functions
\[
\sum_{n=1}^{\infty}H_{n}^{^{-}}t^{n}=\frac{\log\left(  1+t\right)  }%
{1-t}\text{ and }\sum_{n=1}^{\infty}\left(  -1\right)  ^{n}H_{n}t^{n}%
=-\frac{\log\left(  1+t\right)  }{1+t},
\]
we find that%
\[
\int_{0}^{\infty}\frac{x^{s-1}}{1-e^{-x}}\log\left(  1+e^{-x/2}\right)
dx=2^{s-1}\Gamma\left(  s\right)  \left(  \zeta_{H^{^{-}}}\left(  s\right)
+\zeta_{\widetilde{H}}\left(  s\right)  \right)  .
\]
Thus we obtain
\begin{equation}
\frac{1}{2}\sum_{k=1}^{2}\zeta_{H}\left(  s,\frac{k}{2}\right)  =\zeta
_{H}\left(  s\right)  +2^{s-1}\zeta_{H^{^{-}}}\left(  s\right)  +2^{s-1}%
\zeta_{\widetilde{H}}\left(  s\right)  . \label{12}%
\end{equation}
We now use (\ref{12a}), (\ref{12}) and Theorem \ref{mteo2} to conclude that%
\begin{align*}
&  \sum_{m=0}^{\infty}\left(  -1\right)  ^{m}\sum_{k=1}^{2}\frac{1}{2}%
\gamma_{H}\left(  m,\frac{k}{2}\right)  \frac{\left(  s-1\right)  ^{m}}%
{m!}=\sum_{m=0}^{\infty}\left(  -1\right)  ^{m}\gamma_{H}\left(  m\right)
\frac{\left(  s-1\right)  ^{m}}{m!}\\
&  +\sum_{m=0}^{\infty}\left(  \sum_{k=0}^{m}\binom{m}{k}\left(  -1\right)
^{k}\left(  \gamma_{\widetilde{H}}\left(  k\right)  +\gamma_{H^{^{-}}}\left(
k\right)  \right)  \log^{m-k}2+\frac{\log^{m+1}2}{m+1}\right)  \frac{\left(
s-1\right)  ^{m}}{m!}.
\end{align*}
Comparing the coefficients of $\left(  s-1\right)  ^{m}/m!$ we deduce the
desired result.
\end{proof}

Next proposition gives an evaluation formula for the Stietljes
constants\ $\gamma_{H}\left(  n\right)  .$ This formula contributes to the
computations of $\gamma_{H}\left(  n,1/2\right)  $ and the hyperharmonic
Stieltjes constants given in \cite[Theorem 5]{CDK}.

\begin{proposition}
For a non-negative integer $n$ we have%
\begin{align}
\gamma_{H}\left(  n\right)   &  =\left(  -1\right)  ^{n}\zeta^{\left(
n\right)  }\left(  2\right)  +\gamma\left(  n+1\right)  -\left(  -1\right)
^{n}\frac{\psi^{\left(  n+1\right)  }\left(  1\right)  }{n+1}\nonumber\\
&  -\left(  -1\right)  ^{n}\sum_{k=0}^{n}\dbinom{n}{k}\left(  b_{k}%
i_{n-k}+\left(  -1\right)  ^{n-k}\gamma\left(  n-k\right)  \psi^{\left(
k\right)  }\left(  1\right)  \right)  , \label{ev11}%
\end{align}
$\allowbreak$where
\begin{equation}
b_{k}=\left.  \frac{d^{k}}{ds^{k}}\frac{1}{\Gamma\left(  s\right)
}\right\vert _{s=1}\text{ and }i_{k}=\int_{0}^{\infty}\frac{e^{-x}\log\left(
\frac{1-e^{-x}}{x}\right)  \log^{k}x}{1-e^{-x}}dx. \label{bkik}%
\end{equation}

\end{proposition}

\begin{proof}
We focus on the RHS of (\ref{B1}). It is known that the Laurent series
expansions of the functions $\psi\left(  s\right)  \zeta\left(  s\right)
$ and $-\zeta^{\prime}\left(  s\right)  $ in a
neighborhood of $s=1$ are
\begin{align*}
\psi\left(  s\right)  \zeta\left(  s\right)     &=\frac{\psi\left(  1\right)
}{s-1}+\sum_{n=0}^{\infty}\frac{\psi^{\left(  n+1\right)  }\left(  1\right)
}{n+1}\frac{\left(  s-1\right)  ^{n}}{n!}\\
& +\sum_{n=0}^{\infty}\sum_{k=0}^{n}\dbinom{n}{k}\left(  -1\right)
^{n-k}\gamma\left(  n-k\right)  \psi^{\left(  k\right)  }\left(  1\right)
\frac{\left(  s-1\right)  ^{n}}{n!}%
\end{align*}
and%
\[
\zeta^{\prime}\left(  s\right)  =-\frac{1}{\left(  s-1\right)  ^{2}}%
-\sum_{n=0}^{\infty}\left(  -1\right)  ^{n}\gamma\left(  n+1\right)
\frac{\left(  s-1\right)  ^{n}}{n!}.
\]
To expand $\Gamma\left(  1-s\right)  I\left(  s\right)  $ into a series we
first write $I\left(  s\right)  $ as%
\begin{equation}
I\left(  s\right)  =-\frac{\sin\left(  \pi s\right)  }{\pi}\int_{0}^{\infty
}\frac{x^{s-1}e^{-x}}{1-e^{-x}}\log\left(  \frac{1-e^{-x}}{x}\right)  dx,
\label{Is}%
\end{equation}
and $\Gamma\left(  1-s\right)  I\left(  s\right)  $ as%
\[
\Gamma\left(  1-s\right)  I\left(  s\right)  =-\frac{1}{\Gamma\left(
s\right)  }\int_{0}^{\infty}\frac{x^{s-1}e^{-x}}{1-e^{-x}}\log\left(
\frac{1-e^{-x}}{x}\right)  dx.
\]
Since $\frac{1}{\Gamma\left(  s\right)  }$ and the integral on the RHS are
analytic around $s=1$ we can write
\[
\frac{1}{\Gamma\left(  s\right)  }=\sum_{k=0}^{\infty}\frac{b_{k}}{k!}\left(
s-1\right)  ^{k}\text{ and }\int_{0}^{\infty}\frac{x^{s-1}e^{-x}}{1-e^{-x}%
}\log\left(  \frac{1-e^{-x}}{x}\right)  dx=\sum_{k=0}^{\infty}\frac{i_{k}}%
{k!}\left(  s-1\right)  ^{k}.
\]
Hence
\[
\Gamma\left(  1-s\right)  I\left(  s\right)  =-\sum_{n=0}^{\infty}\sum
_{k=0}^{n}\dbinom{n}{k}b_{k}i_{n-k}\frac{\left(  s-1\right)  ^{n}}{n!}.
\]
Thus, we deduce from (\ref{0}) and (\ref{B1}) that%
\begin{align*}
 &\sum_{n=0}^{\infty}\left(  -1\right)  ^{n}\gamma_{H}\left(  n\right)
\frac{\left(  s-1\right)  ^{n}}{n!}\\
&  =\sum_{n=0}^{\infty}\left(  \zeta^{\left(  n\right)  }\left(  2\right)
+\left(  -1\right)  ^{n}\gamma\left(  n+1\right)  -\frac{\psi^{\left(
n+1\right)  }\left(  1\right)  }{n+1}\right)  \frac{\left(  s-1\right)  ^{n}%
}{n!}\\
&  -\sum_{n=0}^{\infty}\sum_{k=0}^{n}\dbinom{n}{k}\left(  b_{k}i_{n-k}+\left(
-1\right)  ^{n-k}\gamma\left(  n-k\right)  \psi^{\left(  k\right)  }\left(
1\right)  \right)  \frac{\left(  s-1\right)  ^{n}}{n!},
\end{align*}
which gives (\ref{ev11}).
\end{proof}

\begin{proposition}
For a non-negative integer $n$ we have%
\[
\gamma_{\widetilde{H}}\left(  n\right)  =\left(  -1\right)  ^{n}\sum_{k=0}%
^{n}\dbinom{n}{k}b_{k}c_{n-k},
\]
where $b_{k}$ is given by (\ref{bkik}) and%
\[
c_{k}=\int_{0}^{\infty}\frac{\log\left(  1+e^{-x}\right)  }{1+e^{-x}}\log
^{k}xdx.
\]

\end{proposition}

\begin{proof}
Similar to (\ref{10}) it can be seen that%
\[
\zeta_{\widetilde{H}}\left(  s\right)  =\frac{1}{\Gamma\left(  s\right)  }%
\int_{0}^{\infty}\frac{x^{s-1}}{1+e^{-x}}\log\left(  1+e^{-x}\right)  dx.
\]
Since $\frac{1}{\Gamma\left(  s\right)  }$ and the integral on the RHS is
analytic around $s=1$, we can write
\[
\frac{1}{\Gamma\left(  s\right)  }=\sum_{k=0}^{\infty}\frac{b_{k}}{k!}\left(
s-1\right)  ^{k}\text{ and }\int_{0}^{\infty}\frac{x^{s-1}}{1+e^{-x}}%
\log\left(  1+e^{-x}\right)  dx=\sum_{k=0}^{\infty}\frac{c_{k}}{k!}\left(
s-1\right)  ^{k}.
\]
Hence,
\[
\zeta_{\widetilde{H}}\left(  s\right)  =\sum_{n=0}^{\infty}\sum_{k=0}%
^{n}\dbinom{n}{k}b_{k}c_{n-k}\frac{\left(  s-1\right)  ^{n}}{n!}.
\]
Combining this with (\ref{12a}) completes the proof.
\end{proof}

\begin{remark}
A few first values of $\gamma_{\widetilde{H}}\left(  n\right)  $\textbf{ }can
be listed as follows.%
\begin{align*}
\gamma_{\widetilde{H}}\left(  0\right)   &  =0.58224052646501\ldots &
\gamma_{\widetilde{H}}\left(  5\right)   &  =0.03566856790329\ldots\\
\gamma_{\widetilde{H}}\left(  1\right)   &  =-\;0.203469139835\ldots &
\gamma_{\widetilde{H}}\left(  6\right)   &  =0.01570725581174\ldots\\
\gamma_{\widetilde{H}}\left(  2\right)   &  =-\;0.069650561786\ldots &
\gamma_{\widetilde{H}}\left(  7\right)   &  =-\;0.012735756486\ldots\\
\gamma_{\widetilde{H}}\left(  3\right)   &  =0.00491231177924\ldots &
\gamma_{\widetilde{H}}\left(  8\right)   &  =-\;0.037846238315\ldots\\
\gamma_{\widetilde{H}}\left(  4\right)   &  =0.03583962465759\ldots &
\gamma_{\widetilde{H}}\left(  9\right)   &  =-\;0.047740308386\ldots
\end{align*}

\end{remark}

We conclude the paper by expressing the integral
\[
J\left(  v\right)  =\frac{1}{\Gamma\left(  v\right)  }\int_{0}^{\infty}%
\frac{x^{v-1}e^{-x}}{1-e^{-x}}\log\left(  \frac{1-e^{-x}}{x}\right)  dx
\]
in terms of known constants.

\begin{proposition}
We have
\[
J\left(  1\right)  =\gamma\left(  1\right)  +\frac{\gamma^{2}\left(  0\right)
-\zeta\left(  2\right)  }{2}%
\]
and for integers $v\geq2$,%
\[
J\left(  v\right)  =-\frac{v}{2}\zeta\left(  v+1\right)  -\psi\left(
v\right)  \zeta\left(  v\right)  -\zeta^{\prime}\left(  v\right)  +\frac{1}%
{2}\sum_{j=1}^{v-2}\zeta\left(  v-j\right)  \zeta\left(  j+1\right)  .
\]

\end{proposition}

\begin{proof}
From (\ref{ev11}) and (\ref{bkik}) we have%
\[
\gamma_{H}\left(  0\right)  =\gamma\left(  1\right)  +\gamma^{2}\left(
0\right)  -\int_{0}^{\infty}\frac{e^{-x}}{1-e^{-x}}\log\left(  \frac{1-e^{-x}%
}{x}\right)  dx.
\]
Combining this with
\[
\gamma_{H}\left(  0\right)  =\frac{\gamma^{2}\left(  0\right)  +\zeta\left(
2\right)  }{2}%
\]
(see \cite[Eq. (6)]{CC} or \cite[Remark 2]{CDK}) gives
\[
J\left(  1\right)  =\gamma\left(  1\right)  +\frac{\gamma^{2}\left(  0\right)
-\zeta\left(  2\right)  }{2},
\]
the first equality.

Since $\zeta_{H}\left(  v\right)  $ is convergent for integers $v\geq2$, we
may write (\ref{B1}) in the following form%
\[
\lim_{s\rightarrow v}\frac{\pi I\left(  s\right)  }{\Gamma\left(  s\right)
\sin\left(  \pi s\right)  }=\zeta_{H}\left(  v\right)  -\zeta\left(
v+1\right)  +\psi\left(  v\right)  \zeta\left(  v\right)  +\zeta^{\prime
}\left(  v\right)  .
\]
Employing (\ref{Is}), the LHS becomes%
\[
\lim_{s\rightarrow v}\frac{\pi I\left(  s\right)  }{\Gamma\left(  s\right)
\sin\left(  \pi s\right)  }=-J\left(  v\right)  .
\]
In the light of the well-known Euler's formula%
\[
\zeta_{H}\left(  v\right)  =\frac{v+2}{2}\zeta\left(  v+1\right)  -\frac{1}%
{2}\sum_{j=1}^{v-2}\zeta\left(  v-j\right)  \zeta\left(  j+1\right)  ,
\]
the proof follows.
\end{proof}

\end{document}